\documentclass{amsart}

\usepackage{amsfonts}
\usepackage{amsmath}
\usepackage{amsthm}
\usepackage{amssymb}
\usepackage[english]{babel}
\usepackage{graphics}
\usepackage{latexsym}
\usepackage{longtable} 
\usepackage{mathrsfs} 
\usepackage{graphicx}
\usepackage{wrapfig} 
\usepackage[pdftex,colorlinks=true]{hyperref}

\newtheorem{theorem}{Theorem}
\newtheorem{corollary}[theorem]{Corollary}
\newtheorem{lemma}[theorem]{Lemma}

\newtheorem{claim}[theorem]{Claim}
\newtheorem{example}[theorem]{Example}

\theoremstyle{definition}
\newtheorem{definition}[theorem]{Definition}

\newcommand{\mC}{\mathcal{C}}
\newcommand{\mL}{\mathcal{L}}

\newcommand{\mF}{\mathcal{F}}

\newcommand{\mD}{\mathcal{D}}

\newcommand{\B}{\mathcal{B}}

\newcommand{\R}{\mathbb{R}}

\newcommand{\N}{\mathbb{N}}

\newcommand{\mB}{\mathbb{B}}

\newcommand{\noi}{\noindent}
\newcommand{\ms}{\medskip}

\newcommand{\al}{\alpha}
\newcommand{\be}{\beta}

\newcommand{\de}{\delta}
\newcommand{\De}{\Delta}
\newcommand{\e}{\varepsilon}

\newcommand{\la}{\lambda}

\newcommand{\Om}{\Omega}
\newcommand{\om}{\omega}

\newcommand{\ze}{\zeta}

\newcommand{\D}{\mathrm{D}} 
 
\newcommand{\larrow}{\longrightarrow}
\newcommand{\ot}{\otimes}

\newcommand{\lmapsto}{\longmapsto}
\newcommand{\ri}{\rightarrow}

\newcommand{\p}{\partial}
\newcommand{\sub}{\subseteq}
\newcommand{\set}{\setminus}
\newcommand{\by}{\times}
\newcommand{\rk}{\mathrm{rk}}

\newcommand{\ess}{\mathrm{ess}}

\newcommand{\dist}{\mathrm{dist}}

\newcommand{\spn}{\mathrm{span}}
\newcommand{\supp}{\mathrm{supp}}

\newcommand{\bt}{\begin{theorem}}\newcommand{\et}{\end{theorem}}
\newcommand{\bd}{\begin{definition}}\newcommand{\ed}{\end{definition}}
\newcommand{\bl}{\begin{lemma}}\newcommand{\el}{\end{lemma}}
\newcommand{\beq}{\begin{equation}}\newcommand{\eeq}{\end{equation}}
\newcommand{\bc}{\begin{claim}}\newcommand{\ec}{\end{claim}}
\newcommand{\bex}{\begin{example}}\newcommand{\eex}{\end{example}}
\newcommand{\bcor}{\begin{corollary}}\newcommand{\ecor}{\end{corollary}}
\newcommand{\bp}{\begin{proof}}\newcommand{\ep}{\end{proof}}

\numberwithin{equation}{section}

\begin{document}

\title[Vectorial Solutions of H-J PDE are rank-one Absolute Minimisers]{Solutions of Vectorial Hamilton-Jacobi Equations are Rank-One Absolute Minimisers in $L^\infty$}


\author{Nikos Katzourakis}

\thanks{\texttt{The author has been partially supported by an EPSRC grant EP/N017412/1.} }

\address{Department of Mathematics and Statistics, University of Reading, Whiteknights, PO Box 220, Reading RG6 6AX, Berkshire, UK}
\email{n.katzourakis@reading.ac.uk}



\date{}



\begin{abstract} Given the supremal functional $E_\infty(u,\Omega')=\ess\,\sup_{\Omega'} H(\cdot,\D u)$ defined on $W^{1,\infty}_{\text{loc}}(\Omega,\mathbb{R}^N)$, $\Omega' \Subset \Omega\subseteq \mathbb{R}^n$, we identify a class of vectorial rank-one Absolute Minimisers by proving a statement slightly stronger than the next claim: vectorial solutions of the Hamilton-Jacobi equation $H(\cdot,\D u)=c$ are rank-one Absolute Minimisers if they are $C^1$. Our minimality notion is a generalisation of the classical $L^\infty$ variational principle of Aronsson to the vector case and emerged in earlier work of the author. The assumptions are minimal, requiring only continuity and rank-one convexity of the level sets. 
\end{abstract}

\maketitle


\section{Introduction} \label{section1}

In this paper we are concerned with the construction of a special class of appropriately defined vectorial minimisers in Calculus of Variations in $L^\infty$, that is for supremal functionals of the form
\beq \label{a1}
\left\{\ \ \
\begin{split}
 & E_\infty (u,\Om')\, :=\, \underset{x\in \Om'}{\ess\,\sup}\, H\big(x,\D u(x)\big), \ \ \ \ \
\\
 & \ \ \ \ \ u  \,\in W^{1,\infty}_{\text{loc}}(\Om,\R^N), \ \  \Om'\Subset \Om.
\end{split}
\right.
\eeq
In the above, $n,N\in \N$, $\Om \sub \R^n$ is an open set and $H : \Om \by \R^{N\by n} \larrow [0,\infty)$ is a continuous function. Calculus of Variations in $L^\infty$ has been pioneered by Aronsson in the 1960s who studied the scalar case $N=1$ quite systematically (\cite{A1}-\cite{A6}). A major difficulty associated to the study of \eqref{a1} is that the standard minimality notion used for the respective more classical integral functional
\beq \label{a2}
 E (u,\Om)\, =\, \int_{\Om} H\big(x,\D u(x)\big)\,dx , 
\eeq
is not appropriate in the $L^\infty$ case due to the lack of ``locality" of \eqref{a1}. The remedy is to require minimality on all subdomains $\Om'\Subset \Om$, a notion now known as Absolute Minimality (hence the emergence of the domain as second argument of the functional). The field has seen an explosion of interest especially after the 1990s when the development of Viscosity solutions (\cite{CIL, C}) allowed the rigorous study of the non-divergence equation arising as the analogue of the Euler-Lagrange for \eqref{a1} (for a pedagogical introduction to the scalar case and numerous references we refer to \cite{K7}). In the special case of $H(x,P)$ being the Euclidean norm $|P|$ on $\R^n$, the respective PDE is the $\infty$-Laplace equation
\beq \label{a3}
\De_\infty u \ :=\ \D u \ot \D u : \D^2 u\ = \sum_{i,j=1}^n\D_i u\, \D_ju\,\D^2_{ij}u\,=\ 0. 
\eeq
Despite the importance for applications and the deep analytical interest of the area, the vectorial case of $N\geq 2$ remained largely unexplored until the early 2010s. In particular, not even the correct form of the respective PDE systems associated to $L^\infty$ variational problem was known. A notable exception is the early vectorial contributions \cite{BJW1, BJW2} wherein (among other deep results) $L^\infty$ versions of lower semi-continuity and quasiconvexity were introduced and studied and the existence of Absolute Minimisers was established in some generality with $H$ depending on $u$ itself but for $\min\{n,N\}=1$. 

The author in a series of recent papers (see \cite{K1}-\cite{K6}, \cite{K8}-\cite{K11} has laid the foundations of the vectorial case and in particular has derived and studied the analogues of \eqref{a3} associated to general $L^\infty$ functionals (see also the joint contributions with Croce, Pisante, Pryer and Abugirda \cite{AK, CKP, KP1, KP2}). In the model case of $H$ being the Euclidean norm on $\R^{N\by n}$ and independent of $x$
\[
H(x,P)\, =\, |P|\, = \left(\sum_{\al=1}^N\sum_{i=1}^n(P_{\al i})^2\right)^{1/2},
\]
the respective equation is called the $\infty$-Laplace PDE system and when applied to smooth maps $u : \Om \sub \R^n \larrow \R^N$ reads
\beq  \label{a4}
\De_\infty u \ :=\ \Big(\D u \ot \D u  + |\D u|^2 [\D u]^\bot \! \ot \mathrm{I} \Big) : \D^2 u\ = \ 0. 
\eeq
Here $[Du(x)]^\bot$ is the orthogonal projection on the orthogonal complement of the range  of gradient matrix $\D u(x) \in \R^{N\by n}$: 
\[
[\D u]^\bot :=\, \mathrm{Proj}_{(R(\D u))^\bot}.
\]
In index form,  \eqref{a4} reads
\[  
\ \ \ \ \sum_{\be=1}^N \sum_{i,j=1}^n \Big(\D_i u_\al  \D_j u_\be  \, +\, |\D u|^2 [\D u]_{\al \be}^\bot\,  \de_{ij}\Big) \, \D^2_{ij} u_\be\ = \ 0, \ \ \ \ \al=1,...,N. 
\]

In the full vector case of \eqref{a4}, even more intriguing phenomena occur since a further difficulty which is not present in the scalar case is that the coefficient involving $[\D u]^\bot$ is discontinuous even for $C^\infty$ maps; for instance, $u(x,y) = e^{ix}-e^{iy}$ is a smooth $2\by 2$ $\infty$-Harmonic map near the origin and the rank of gradient is $1$ on the diagonal, but it is $2$ otherwise. The emergence of discontinuities is a genuine vectorial phenomenon which does not arise if $\min\{n,N\}=1$ (see \cite{K1,K3,K4}). For $N=1$ the scalar version \eqref{a3} has continuous coefficients, whilst for $n=1$ \eqref{a4} reduces to 
\[
\De_\infty u \, =\, (u' \ot u') u''\, +\, |u'|^2\Big(I - \frac{u'}{|u'|}\ot \frac{u'}{|u'|}\Big)u''\ =\ |u'|^2u''.
\]
A problem associated to the discontinuities is that \emph{Aronsson's notion of Absolute Minimisers is not appropriate in the vectorial case of rank $\rk(\D u)\geq 2$}. Actually, by the perpendicularity of $\D u$ and $[\D u]^\bot$, $\De_\infty u=0$ actually consists of $2$ independent systems and each one is characterised in terms of the $L^\infty$ norm of the gradient via different sets of variations. In \cite{K2} we proved the following variational characterisation in the class of classical solutions. A $C^2$ map $u :\Om\sub \R^n \larrow \R^N$ is a solution to
\beq \label{a5}
\D u \ot \D u :\D^2u\, =\,0
\eeq
if and only if it is a \emph{Rank-One Absolute Minimiser} on $\Om$, namely when for all  $D \Subset \Om$, all scalar functions $g \in C^1_0(D)$ vanishing on $\p D$ and all directions $\xi \in \R^N$, $u$ is a minimiser on $D$ with respect to variations  of the form $u+\xi g$ (Figure 1):
\beq \label{a6}
\|\D u \|_{L^\infty(D)}\ \leq \ \big\|\D u+ \xi \ot \D g\big\|_{L^\infty(D)}.
\eeq
\[
\underset{\text{Figure 1.}}{\includegraphics[scale=0.18]{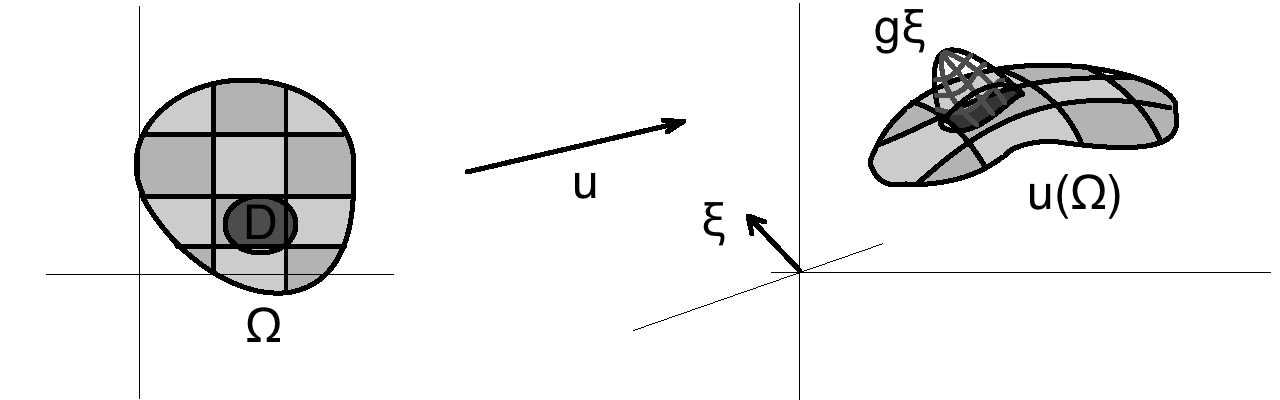}}
\]
Further, if $\rk(\D u)\equiv$ const., $u$ is a solution to 
\beq  \label{a7}
|\D u|^2 [\D u]^\bot \De u \,=\, 0 
\eeq
if and only if \emph{$u(\Om)$ has $\infty$-Minimal Area}, namely when for all $D \Subset\Om$, all scalar functions $h  \in C^1(\overline{D})$ (not vanishing on $\p D$) and all vector fields $\nu \in C^1(D,\R^N)$ which are normal to $u(\Om)$, $u$ is a minimiser on $D$ with respect to normal free variations of the form $u+h\nu$ (Figure 2):
\beq \label{a8}
\|\D u \|_{L^\infty(D)}\ \leq \ \big\|\D u+\D(h\nu)\big\|_{L^\infty(D)}.
\eeq
\[
\underset{\text{Figure 2.}}{\includegraphics[scale=0.18]{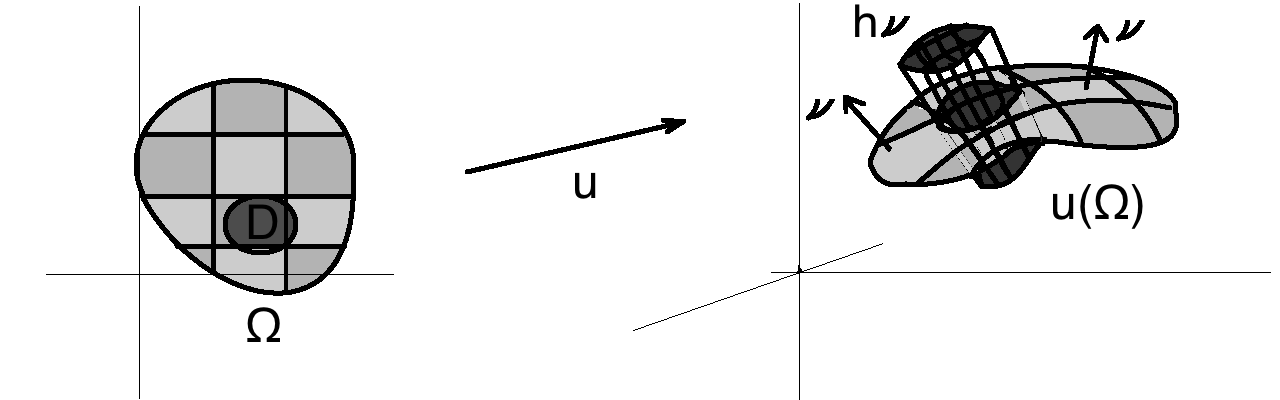}}
\]
We called a map $\infty$-Minimal with respect to functional $\|\D(\cdot)\|_{L^\infty(\cdot)}$ when it is a Rank-One Absolute Minimiser on $\Om$ and $u(\Om)$ has $\infty$-Minimal Area.

Perhaps the greatest difficulty associated to \eqref{a1} and \eqref{a4} is how to define and study generalised solutions since for the highly nonlinear non-divergence model system \eqref{a4} all standard arguments based on the maximum principle or on integration by parts seem to fail. In the very recent work \cite{K8} the author proposed the theory of so-called \emph{$\mD$-solutions} which applies to general fully nonlinear PDE system of any order
\[
\mF\Big(\cdot,u,\D u,...,\D^pu\Big)\,=\,0
\]
and allows for merely measurable solutions $u$ to be rigorously interpreted and studied. This notion is duality-free and is based on the probabilistic representation of derivatives which do not exist classically. $\mD$-solutions have already borne substantial fruit and in \cite{K8}-\cite{K12} we have derived several existence-uniqueness, variational and regularity results. In particular, in \cite{K10} we have obtained a variational characterisation of \eqref{a4} in the setting of general $\mD$-solutions for $W^{1,\infty}_{\text{loc}}(\Om,\R^N)$ appropriately defined minimisers which is relevant to \eqref{a5}-\eqref{a8} but different.

In this paper we consider the obvious generalisation of the rank-one minimality notion of \eqref{a6} adapted to the functional \eqref{a1}. To this end, we identify a large class of rank-one Absolute Minimisers: for any $c\geq0$, every solution $u : \Om\sub \R^n \larrow \R^N$ to the vectorial Hamilton-Jacobi equation
\beq
H\big(x,\D u(x)\big)\,=\,c,\ \ \ \ x\in \Om,
\eeq
actually is a rank-one absolute minimiser. Namely, for any $\Om'\Subset \Om$, any $\phi \in W_0^{1,\infty}(\Om')$ and any $\xi \in \R^N$, we have
\[
\underset{x\in \Om'}{\ess\,\sup}\, H\big(x,\D u(x)\big)\, \leq \, \underset{x\in \Om'}{\ess\,\sup}\, H\Big(x,\D u(x)+\xi \ot \D \phi(x)\Big).
\]
For the above implication to be true we need the solutions to be in $C^1(\Om,\R^N)$ and not just in $W^{1,\infty}_{\text{loc}}(\Om,\R^N)$. This is not a technical difficulty: it is well known even in the scalar case that if we allow only for one non-differentiability point, strong solutions of the Eikonal equation $|\D u|=1$ are not absolutely minimising for the $L^\infty$ norm of the gradient (e.g.\ the cone function $x\mapsto |x|$). However, due to regularity results which are available in the scalar case, it suffices to assume everywhere differentiability (see \cite{CEG, CC}). 

Our only hypothesis imposed on $H$ is that for any $x\in \Om$ the partial function $H(x,\cdot) : \R^{N\by n} \larrow \R$ is \emph{rank-one level-convex}. This means that for any $t\geq 0$, the sublevel sets $\big\{H(x,\cdot)\leq t\big\} $ are rank-one convex sets in $\R^{N\by n}$. A set $\mC\sub \R^{N\by n}$ is called rank-one convex when for any matrices $A,B \in \mC$ with $\rk(A-B)\leq1$, the convex combination $\la A +(1-\la)B $  is in $\mC$ for any $0\leq \la\leq 1$ (cf.\ \cite[Def.\ 3.7]{BJW2}). An equivalent way to phrase the rank-one level-convexity of $H(x,\cdot)$ is via the inequality
\[
H\Big(x,\la A +(1-\la)B\Big)\, \leq\, \max\big\{ H(x,A),H(x,B)\big\},
\]
when 
\[
\text{$x\in \Om$, $A,B\in \R^{N\by n}$, $\rk(A-B)\leq1$ and $0\leq \la\leq 1$.}
\]
This convexity assumption is substantially weaker than the $L^\infty$ versions of quasiconvexity which  we call ``BJW-quasiconvexity", named after Barron-Jensen-Wang who introduced it in \cite{BJW2}. 

We note that Hamilton-Jacobi equations are very important for $L^\infty$ variational problems and their equations. In the scalar case, $C^1$ solutions are viscosity solutions to the respective second order single equations in $L^\infty$ (see e.g.\ \cite{K7}). Heuristically, for the case of the $\infty$-Laplace equation this can be seen by rewriting \eqref{a3} as $\D u \D\big(\frac{1}{2}|\D u|^2\big)=0$ and this reveals that solutions of $|\D u|=c$ are $\infty$-Harmonic. In the vectorial case, Hamilton-Jacobi equations give rise to certain first-order differential inclusions of the form
\[
\D u(x) \, \in \, \mathcal{K} \sub \R^{N\by n},\ \ \ x\in \Om,
\]
for which the Dacorogna-Marcellini Baire category (the analytic counterpart of Gromov's Convex Integration) method can be utilised to establish existence of $\mD$-solutions to the $L^\infty$ systems of PDE with extra geometric properties (see \cite{K8} and \cite{DM, D}). 

The main result of the present paper is the following:

\begin{theorem} \label{theorem1} Let $\Om\sub \R^n$ be an open set, $n,N\in \N$ and $H : \Om \by \R^{N\by n} \larrow [0,\infty)$ a continuous function, such that for all $x\in \Om$, $P\lmapsto H(x,P)$ is rank-one level-convex, that is
\[
\big\{H(x,\cdot)\leq t\big\} \text{ is a rank-one convex in $\R^{N\by n}$, for all }t\geq 0,\, x\in \Om.
\]
Let $u\in C^1(\Om,\R^N)$ be a solution to the vectorial Hamilton-Jacobi PDE
\[
H(\cdot,\D u)\, =\, c\ \ \text{ on }\Om,
\]
for some $c\geq 0$. Then, $u$ is a rank-one Absolute Minimiser of the functional
\[
E_\infty(u,\Om')\, =\, \underset{x\in \Om'}{\ess\,\sup}\, H\big(x,\D u(x)\big), \ \ \ \Om'\Subset \Om,\ u\in W^{1,\infty}_{\text{loc}}(\Om,\R^N).
\]
In addition, the following marginally stronger result holds true: for any $\Om'\Subset \Om$, any $\phi \in W_0^{1,\infty}(\Om')$ and any $\xi \in \R^N$, we have
\[
E_\infty(u,\Om')\, \leq \,  \inf_{\mB\in \B(\phi,\Om')} \, E_\infty\big(u+\xi \phi,\mB\big)
\]
where $\B(\phi,\Om')$ is the set of open balls centred at local extrema (maxima or minima) of $\phi$ inside $\Om'$:
\beq \label{1.5}
\B(\phi,\Om') \,:=\, \Big\{\mB_\rho(x) \Subset \Om'\, \Big|\ x \text{ is a point of local extremum of }\phi \text{ in }\Om'\Big\}.
\eeq
\end{theorem}

An immediate consequence of Theorem \ref{theorem1} is the following result:

\begin{corollary} In the setting of Theorem \ref{theorem1}, we additionally have
\[
\underset{ \Om'}{\ess\,\sup}\, H\big(\cdot,\D u\big)\, \leq \,  \underset{\rho\ri 0}{\lim} \bigg( \underset{ \mB_\rho(x)}{\ess\,\sup}\, H\big(\cdot,\D u +\xi \ot \D \phi\big)\bigg),
\]
for any $x\in \Om'$ at which $\phi$ achieves a local maximum or a local minimum. \end{corollary}

The quantity of the right hand side above is known as the local functional at $x$ and in the scalar case has been used as a substitute of the pointwise values due to its upper semi-continuity regularity properties.

\section{The proof of Theorem \ref{theorem1}.}

Let $H,c,u,\Om$ be as in the statement and fix $\Om'\Subset \Om$ and a unit vector $\xi \in \R^N$. We introduce the following notation for the projections on $\spn[\xi]$ and the orthogonal hyperplane $\spn[\xi]^\bot$:
\beq
\left\{\ \ 
\begin{split}
[\xi]^\top\, &:=\, \xi \ot \xi ,\\
[\xi]^\bot\, &:=\, I-\xi \ot \xi.  \ \ \ \ 
\end{split} 
\right.
\eeq
Let $\psi \in W^{1,\infty}_u(\Om',\R^N):=u+W^{1,\infty}_0(\Om',\R^N)$ such that $[\xi]^\bot\big( \psi - u\big) \equiv  0$, that is the projections of $\psi$ and $u$ on the hyperplane $\spn[\xi]^\bot\sub \R^N$ coincide. Then, $\xi\cdot(\psi-u) \in W^{1,\infty}_0(\Om')$ and because the scalar function $\xi\cdot(\psi-u )$ vanishes on $\p\Om'$, there exist at least one local extremum of $\xi\cdot(\psi-u )$ in $\Om'$, whence the set 
\[
\B\big(\xi \cdot(\psi-u),\Om'\big)
\]
given by \eqref{1.5} is non-empty. Fix a ball $\mB_\rho(x)\Subset \Om'$ centred at such an extremal point of $\xi\cdot(\psi-u)$. 

We illustrate the idea by assuming first in addition that $\psi \in W^{1,\infty}_u(\Om',\R^N)\cap C^1(\Om',\R^N)$. In this case, the point $x$ is a critical point of $\xi\cdot(\psi-u)$ and we have $\D \big(\xi\cdot(\psi-u)\big)(x)=0$. Hence,
\[
\begin{split}
\D(\psi-u)(x)\,&=\, [\xi]^\top \D(\psi-u)(x)\, +\, [\xi]^\bot \D(\psi-u)(x)
\\
&=\, \xi \ot \D \big(\xi\cdot(\psi-u)\big)(x)\, +\, \D \big([\xi]^\bot (\psi-u)\big)(x)
\\
&=\,0
\end{split}
\]
because $[\xi]^\bot\psi \equiv [\xi]^\bot u$ on $\Om'$. Thus,
\[
\begin{split}
E_\infty(u,\Om')\, =\, c\, &=\, H(x,\D u(x))\\
&=\, H\Big(x,\, \xi \ot \D (\xi\cdot u) (x)\, +\, [\xi]^\bot \D u(x)\Big)
\end{split}
\]
and hence
\beq \label{1.6}
\begin{split}
E_\infty(u,\Om')\, &=\, H\Big(x,\, \xi \ot \D (\xi\cdot \psi )(x)\, +\, [\xi]^\bot \D \psi (x)\Big)\\
&=\, H(x,\D \psi(x))\\
&\leq\, \underset{y\in \mB_\rho(x)}{\ess\,\sup} \, H(y,\D \psi(y))\\
&= \, E_\infty\big(\psi,\mB_\rho(x)\big)
\end{split}
\eeq
for any $\mB_\rho(x) \sub \B\big(\xi \cdot(\psi-u),\Om'\big)$, whence the conclusion ensues.

Now we return to the general case of $\psi \in W^{1,\infty}_u(\Om',\R^N)$. We extend $\psi$ by $u$ on $\Om\set \Om'$ and consider the sets
\beq \label{1.7}
\Om_k\, :=\, 
\left\{
\begin{array}{ll}
\left\{x\in \Om'\,:\, \dist(x,\p\Om')\,>\, \dfrac{d_0}{k}\right\}, & k\in \N, \ms
\\
\ \emptyset, & k=0,
\end{array}
\right.
\eeq
where $d_0>0$ is a constant small enough so that $\Om_1 \neq \emptyset$. We set
\beq \label{1.8}
V_k\,:=\, \Om_k\set \overline{\Om_{k-1}}, \ \ \ \ k\in \N
\eeq
and consider a partition of unity $(\ze_k)_{k=1}^\infty \sub C^\infty_c(\Om')$ over $\Om'$ so that
\beq \label{1.9}
\left\{\ \ 
\begin{split}
\ze_k \,\geq\, 0\text{ on }\Om', \ \ \ \sum_{k=1}^\infty \ze_k \,\equiv\, 1\text{ on }\Om', 
\\ 
 \ze_k \,>\, 0\text{ on }V_k, \hspace{50pt}
 \\
  \  \supp(\ze_k)\,\sub\,V_{k-1}\cup V_k \cup V_{k+1}. \ \ \ \
\end{split}
\right.
\eeq
(Such a partition of unity can be easily constructed explicitly by mollifying the characteristic functions $\chi_{V_k}$ and rescaling them appropriately.) Let $(\eta^\e)_{\e>0}$ be the standard mollifier (as e.g.\ in \cite[Appendix C.5]{E}) and set 
\beq \label{1.10}
\psi^\e\, :=\, \xi \ot \left(\sum_{k=1}^\infty\, \ze_k \big((\xi \cdot \psi)*\eta^{\e/k}\big)\right) +\, [\xi]^\bot \psi , \ \ \ 0<\e<d_0.
\eeq
We claim that $\psi^\e \in W^{1,\infty}_u(\Om',\R^N) \cap C^1(\Om',\R^N)$ and $\psi^\e \larrow \psi$ in $C(\overline{\Om'},\R^N)$ as $\e\ri0$. To verify these claims, fix $l\in \N$. If $l \geq 2$, then by \eqref{1.8}-\eqref{1.9} and because $|\xi|=1$, we have
\beq \label{1.11}
\begin{split}
\|\psi^\e -\, \psi\|_{C(\overline{V_l})}\, &= \, \sup_{V_l} \left|\xi \cdot\left(\sum_{k=1}^\infty\, \ze_k \big(\psi*\eta^{\e/k}) - \, \psi \sum_{k=1}^\infty\, \ze_k \right)\right|
\\
&\leq \, \sup_{V_l} \sum_{k=1}^\infty\, \ze_k \, \big|  \psi*\eta^{\e/k} -\, \psi \big|
\\
&=\, \sup_{V_l} \sum_{k=l-1}^{l+1}\, \ze_k \, \big|  \psi*\eta^{\e/k} - \, \psi \big|
\\
&\leq \, 3 \max_{k=l-1,l,l+1} \, \big\|  \psi*\eta^{\e/k} -\, \psi \big\|_{C(\overline{\Om'})},
\end{split}
\eeq
whilst, for $l=1$ we similarly have
\beq \label{1.12}
\begin{split}
\|\psi^\e - \, \psi\|_{C(\overline{V_1})}\, \leq \, 2 \max_{k=1,2} \, \big\|  \psi*\eta^{\e/k} - \, \psi \big\|_{C(\overline{\Om'})}.
\end{split}
\eeq
By the standard properties of mollifiers, we have that the function
\beq \label{1.13}
\om(t)\, :=\, \sup_{0<\tau<t}\, \big\|  \psi*\eta^{\tau} - \, \psi \big\|_{C(\overline{\Om'})}, \ \ \ \ 0<t<d_0,
\eeq
is an increasing continuous modulus of continuity with $\om(0^+)=0$. By \eqref{1.11}-\eqref{1.13}, we have that
\beq  \label{1.14}
\|\psi^\e -\, \psi\|_{C(\overline{V_l})}\, \, \leq\, 
\left\{
\begin{array}{ll}
3\, \om\Big(\dfrac{\e}{l-1} \Big), & \ l \geq 2,
\\
2 \, \om(\e), & \ l=1.
\end{array}
\right.
\eeq
Since the $C^1$ regularity of $\psi^\e$ is obvious (because $u$ by assumption is such and $[\xi]^\bot\psi \equiv [\xi]^\bot u$), the claim has been established. 

Note now that since $\psi-u \in W_0^{1,\infty}(\Om',\R^N)$, the set $\B\big(\xi\cdot(\psi-u),\Om'\big)$ given by \eqref{1.5} is non-empty because the scalar function $\xi\cdot(\psi-u)$ which vanishes on $\p\Om'$ necessarily attains an interior extremum. Fix a ball 
\[
\mB_{\rho}(x_0) \, \in \, \B\big(\xi\cdot(\psi-u),\Om'\big).
\]
Since $\psi^\e-u \larrow \psi-u$ in $C(\overline{\Om'},\R^N)$ as $\e\ri0$, by a standard stability argument of maxima/minima of scalar-valued function under uniform convergence (see e.g.\ \cite{K8}), there exists a local extremum $x_\e \in \Om'$ of $\xi \cdot( \psi^\e-u)$ such that $x_\e \larrow x_0$ as $\e \ri 0$. By the differentiability of $u$ and by choosing $\e$ small enough, we may arrange 
\beq  \label{1.15}
\D\big(\xi \cdot (\psi^\e -\, u)\big)(x_\e)\, =\, 0, \ \ \ |x_\e -\,x|\, <\, \frac{\rho}{2}.
\eeq
Hence,
\[
\D(\psi^\e -\, u)(x_\e)\, =\, \xi \ot \D\big(\xi \cdot (\psi^\e -\, u)\big)(x_\e)\, + \,\D\big([\xi]^\bot (\psi\, -\, u)\big)(x_\e)=\, 0.
\]
Then, by arguing as in \eqref{1.6}, we have
\beq  \label{1.16}
\underset{\Om'}{\ess\,\sup} \, H(\cdot,\D u)\, \leq\, \underset{\mB_{\rho/2}(x_0)}{\ess\,\sup} \, H(\cdot,\D \psi^\e).
\eeq
Since
\[
\D \psi^\e\, =\, \xi \ot\left[\sum_{k=1}^\infty\, \D \ze_k \big((\xi \cdot \psi)*\eta^{\e/k}\big) \,+\, \sum_{k=1}^\infty\, \ze_k \Big(\big(\D (\xi \cdot \psi)\big)*\eta^{\e/k}\Big)\right] +\, [\xi]^\bot\D \psi
\]
our continuity assumption and the $W^{1,\infty}$ regularity of $\psi$ imply that there exists a positive increasing modulus of continuity $\om_1$ with $\om_1(0^+)=0$ such that on the ball $\mB_{\rho/2}(x_0)$ we have
\beq \label{1.18}
\begin{split}
H(\cdot,\D \psi^\e) \, &=\, H \Bigg(\cdot\, , \,  \xi \ot \Bigg[ \sum_{k=1}^\infty\, \ze_k \Big(\big(\D (\xi \cdot \psi)\big)*\eta^{\e/k}\Big)  \\
& \hspace{60pt} +\, \sum_{k=1}^\infty\, \D \ze_k \big((\xi \cdot \psi)*\eta^{\e/k}\big) \Bigg] +\, [\xi]^\bot\D \psi \Bigg)
\\
& \leq \, H \left(\cdot\, ,\, \xi \ot \Bigg[\sum_{k=1}^\infty\, \ze_k \big(\D(\xi \cdot \psi)*\eta^{\e/k}\big)\Bigg]  +\, [\xi]^\bot\D \psi \right)\\
& \ \ \ \ +\,  \om_1\left(\left|\sum_{k=1}^\infty\, \D \ze_k \big((\xi \cdot \psi)*\eta^{\e/k}) \right|\right)
\\
&=:\, A\,+\, B.
\end{split}
\eeq
We are planning to show that $B$ tends to zero with $\e$. By further restricting $\e<\rho/2$, we may arrange 
\beq  \label{1.19}
\bigcup_{x\in \mB_{\rho/2}(x_0)}\mB_\e(x)\, \sub\, \mB_{\rho}(x_0)
\eeq
and by \eqref{1.8}-\eqref{1.9}, there exists $K(\rho) \in \N$ such that
\beq  \label{1.20}
\mB_{\rho}(x_0) \,\sub \, \bigcup_{k=1,...,K(\rho)} \overline{V_k}.
\eeq
This implies that for any $x\in \mB_{\rho}(x_0)$,
\beq \label{1.21}
\sum_1^\infty\ze_k(x) \, =\, \sum_1^{K(\rho)+1}\ze_k(x)=1
\eeq
forming a convex combination. We now recall for immediate use right below the following Jensen-like inequality for level-convex functions (see e.g.\ \cite{BJW1, BJW2}): for any probability measure $\mu$ on an open set $U \sub \R^n$ and any $\mu$-measurable function $f :U\sub \R^n \larrow [0,\infty)$, we have
\beq  \label{1.22}
\Phi\left( \int_U f(x)\,d\mu(x) \right) \,\leq\, \mu-\underset{x\in U}{\ess\,\sup} \, \Phi\big(f(x)\big),
\eeq
when $\Phi : \R^n \larrow \R$ is any continuous level-convex function. Further, by our rank-one level-convexity assumption on $H$ and if $\psi$ is as above, for any $x\in \Om$ and $\xi \in \R^N$ with $|\xi|=1$, the function
\beq \label{Psi}
\Psi(p)\,:=\, H\Big(x\, ,\, \xi \ot p \, +\, [\xi]^\bot\D \psi(x) \Big), \ \ \ \ p\in \R^n,
\eeq
is level-convex. Indeed, given $p,q \in \R^n$ and $t\geq 0$ with $\Psi(p),\Psi(q)\leq t$, we set 
\[
\left\{\ \ 
\begin{split}
P\,&:=\, \xi \ot p \, +\, [\xi]^\bot\D \psi(x),\\
Q\,&:=\, \xi \ot q \, +\, [\xi]^\bot\D \psi(x). \ \ \ \ 
\end{split}
\right.
\]
Then, $P-Q=\xi \ot(p-q)$ and hence $\rk(P-Q)\leq 1$. Moreover, $H(x,P)=\Psi(p) \leq t$ and $H(x,Q)=\Psi(q) \leq t$ which gives
\[
\Psi\big(\la p \,+\, (1-\la)q\big)\, =\, H\Big(x,\la P \,+\, (1-\la)Q\Big)\, \leq\, t 
\]
for any $\la \in [0,1]$, as desired.

Now, by using \eqref{1.8}-\eqref{1.9}, \eqref{1.19}-\eqref{1.21} and the level-convexity of the function $\Psi$ of \eqref{Psi}, for any $x\in \mB_{\rho/2}(x_0)$ we have the estimate
\beq \label{1.23}
\begin{split}
A(x)\, &=\, H \left(x\, ,\, \xi \ot \Bigg[\sum_{k=1}^{K(\rho)+1}\, \ze_k(x) \Big(\big(\D (\xi\cdot\psi)\big)*\eta^{\e/k}\Big)(x)\Bigg] +\, [\xi]^\bot\D \psi(x) \right)
\\
&=\, \Psi \left( \sum_{k=1}^{K(\rho)+1}\, \ze_k(x) \Big(\big(\D (\xi\cdot\psi)\big)*\eta^{\e/k}\Big)(x) \right)
\\
& \leq\, \max_{k=1,...,K(\rho)+1}\, \Psi \left( \big( \big(\D (\xi\cdot\psi)\big)*\eta^{\e/k}\big)(x) \right)
\\
& =\, \max_{k=1,...,K(\rho)+1}\, \Psi \left( \int_{\mB_{\e/k}(x)}\D (\xi\cdot\psi)(y)\,\eta^{\e/k}(|x-y|)\,dy \right).
\end{split}
\eeq
Since for any $x$ and $\e,k$, the map
\[
\mu\, :=\, \eta^{\e/k}(|x - \cdot|)\, \mL^n
\]
is a probability measure on the ball $\mB_{\e/k}(x)$ which is absolutely continuous with respect to the Lebesgue measure $ \mL^n$, in view of \eqref{1.22}, \eqref{1.23} gives
\beq \label{1.23a}
\begin{split}
A(x)\, 
& \leq\, \max_{k=1,...,K(\rho)+1}\left( \underset{y\in \mB_{\e/k}(x)}{\ess\,\sup} \, \Psi \Big( \D (\xi \cdot \psi) (y) \Big)\right)
\\
& =\, \max_{k=1,...,K(\rho)+1}\left( \underset{y\in \mB_{\e/k}(x)}{\ess\,\sup} \, H \Big(x\, , \xi \ot \D (\xi \cdot \psi) (y)  +\, [\xi]^\bot\D \psi(x) \Big)\right)
\\
& \leq\,\underset{y\in \mB_{\e}(x)}{\ess\,\sup}  \ H \Big(x\, , \, \xi \ot \D (\xi \cdot \psi) (y)  +\, [\xi]^\bot\D \psi(x) \Big) .
\end{split}
\eeq
By the continuity of $H$ and $\D u$, there is a positive increasing modulus of continuity $\om_2$ with $\om_2(0^+)=0$ such that
\[
\left\{\ \ \ 
\begin{split}
\Big| H(x,P)\,-\, H(y,Q)\Big|\,&\leq\, \om_2\Big(|x-y|\,+\, |P-Q|\Big),
\\
\big|\D u(x)\,-\,\D u(y)\big|\,&\leq\, \om_2\big(|x-y|\big), \ \ \ \ \ \ 
\end{split}
\right.
\]
for all $x,y \in \mB_{\rho}(x_0)$ and $|P|,|Q|\leq\|\D \psi\|_{L^\infty(\Om')}+1$. By using that $[\xi]^\bot \psi \equiv [\xi]^\bot u$ on $\Om'$, \eqref{1.23a} and the above give
\beq \label{1.24}
\begin{split}
A(x)\, & \leq\,\underset{y\in \mB_{\e}(x)}{\ess\,\sup}  \ H \Big(x , \, [\xi]^\top \D \psi (y)  +\, [\xi]^\bot\D \psi(x) \Big) 
\\
 & \leq\,\underset{y\in \mB_{\e}(x)}{\ess\,\sup} \ \Bigg\{ H \Big(y, \, [\xi]^\top \D \psi (y)  +\, [\xi]^\bot\D \psi(y) \Big) 
 \\
 &\hspace{50pt} +\, \om_2\Big(|x-y|\,+\, \Big|[\xi]^\bot\big(\D \psi(y)-\D \psi(x)\big)\Big|\Big)\Bigg\}
 \\
 &= \,\underset{y\in \mB_{\e}(x)}{\ess\,\sup} \ \Bigg\{ H \big(y, \,  \D \psi(y) \big) \, + \om_2\Big(|x-y|\,+\, \Big|[\xi]^\bot\big(\D u(y)-\D u(x)\big)\Big|\Big)\Bigg\}.
\end{split}
\eeq
By \eqref{1.19}, \eqref{1.24} gives
\beq \label{1.25}
\begin{split}
A(x)\, &\leq \,\underset{y\in \mB_{\e}(x)}{\ess\,\sup} \, H \big(y, \,  \D \psi(y) \big) \, + \underset{y\in \mB_{\e}(x)}{\sup} \,\om_2\Big(|x-y|\,+\, \big|\D u(y)-\D u(x)\big|\Big) 
\\
&\leq\, \underset{y\in \mB_{\rho}(x_0)}{\ess\,\sup} \, H \big(y\, , \D \psi (y) \big) \,+\, \om_2\big(\e\,+\, \om_2(\e) \big),
\end{split}
\eeq
for any $x\in \mB_{\rho/2}(x_0)$. We now estimate the term $B$ of \eqref{1.18} on $\mB_{\rho/2}(x_0)$ as above and by using \eqref{1.13}:
\[
\begin{split}
B\, 
&=\, \om_1\left(\left|\xi \cdot \left[\sum_{k=1}^\infty\, \D \ze_k \big(\psi*\eta^{\e/k}) \right] \right|\right)
\\
&\leq\, \om_1\left(\left|\sum_{k=1}^\infty\, \D \ze_k \big(\psi*\eta^{\e/k}) \right|\right)
\\
&\leq\,  \om_1\left(\left|\sum_{k=1}^\infty\, \D \ze_k \,\psi \right| \, +\, \sum_{k=1}^{K(\rho)+1} |\D \ze_k| \Big| \psi*\eta^{\e/k} -\, \psi\Big|\right)
\end{split}
\]
and since $\sum_{k=1}^\infty\, \D \ze_k \equiv 0$, we get
\beq \label{1.26}
\begin{split}
B\, 
&\leq\,  \om_1\left( \sum_{k=1}^{K(\rho)+1} |\D \ze_k| \Big| \psi*\eta^{\e/k} -\, \psi\Big|\right)
\\
&\leq \,  \om_1\left( C(\rho) \max_{k=1,...,K(\rho)+1} \Big\| \psi*\eta^{\e/k} -\, \psi\Big\|_{C(\overline{\Om'})}\right)
\\
& \leq\, \om_1\big( C(\rho) \, \om(\e)\big),
\end{split}
\eeq
on $\mB_{\rho/2}(x_0)$. By putting together \eqref{1.16}, \eqref{1.18}, \eqref{1.25} and \eqref{1.26}, we have
\[
\begin{split}
\underset{\Om'}{\ess\,\sup} \, H(\cdot,\D u)\, &\leq\, \underset{\mB_{\rho}(x_0)}{\ess\,\sup} \, H \big(\cdot\, , \D \psi \big) 
\\
&\ \ \ \ +\, \om_1\big( C(\rho) \, \om(\e)\big) \,+\, \om_2\big(\e\,+\, \om_2(\e) \big)
\end{split}
\]
and by letting $\e\ri 0$, the conclusion follows. 
\qed

\ms
\ms

\noi {\bf Acknowledgement.} The author has been partially financially supported by the EPSRC grant EP/N017412/1.

\ms
\ms

\bibliographystyle{amsplain}

\end{document}